\documentclass[12pt]{amsart}
\usepackage{amssymb}

\newtheorem{thm}{Theorem}[section]
\newtheorem{cor}[thm]{Corollary}
\newtheorem{lem}[thm]{Lemma}

\theoremstyle{definition}

\theoremstyle{remark}
\newtheorem{rem}[thm]{Remark}
\numberwithin{equation}{section}

\newcommand{\eps}{\varepsilon}

\newcommand{\N}{\mathbb{N}}

\newcommand{\C}{\mathbb{C}}
\newcommand{\ov}{\overline}
\newcommand{\dis}{\displaystyle}
\newcommand{\Notequiv}{/\kern-.6em\hbox{$\equiv$} }
\newcommand{\ole}{\preccurlyeq}
\newcommand{\oge}{\succcurlyeq}

\begin{document}

\title[Bieberbach polynomials in domains with interior cusps]
{Convergence of Bieberbach polynomials in domains with interior cusps}%
\author{V. V. Andrievskii and I. E. Pritsker}%
\address{Institute of Biomathematics and Biometry,
GSF-National Research Center for Environment and Health,
D-85764 Neuherberg, Germany}
\email{mgk002@ku-eichstaett.de}%

\address{Department of Mathematics, 401 Mathematical Sciences, Oklahoma State
University, Stillwater, OK 74078-1058, U.S.A.}%
\email{igor@math.okstate.edu}

\thanks{Research of both authors was supported in part by the National Science
Foundation grant DMS-9707359. Research of the second author was
also supported in part by the National Science Foundation
grant DMS-9970659.}%
\subjclass{30E10, 41A10, 30C40}%
\keywords{Bieberbach polynomials, conformal mapping, uniform convergence}%

\begin{abstract}

We extend the results on the uniform convergence of Bieberbach
polynomials to domains with certain interior zero angles (outward
pointing cusps), and show that they play a special role in the
problem. Namely, we construct a Keldysh-type example on the
divergence of Bieberbach polynomials at an outward pointing cusp
and discuss the {\em critical order of tangency} at this interior
zero angle, separating the convergent behavior of Bieberbach
polynomials from the divergent one for sufficiently thin cusps.

\end{abstract}
\maketitle
\section{Introduction}

Let $G$ be a bounded Jordan domain, $z_0 \in G$.  Define the Bergman space $L_2
(G)$ as the space of square integrable analytic functions with norm $$\| f
\|_{2} := \left( \int \!\!\! \int_G  |f(z)|^2 dx dy \right)^{1/2}.$$ We also
use the uniform norm on $\ov G$ in the sequel: $$ \|f\|_{\infty} := \sup_{z \in
G} |f(z)|.$$ The Bieberbach polynomial $B_n (z), \deg B_n \le n$, is the
solution of the following extremal problem \cite{Bi}:
\begin{equation} \label{5.1}
\| B_n' \|_{2} = \inf_{P_n \in {\mathcal P}_n ({\C})}
\{ \| P'_n \|_{2}\ :\ P_n ( z_0 )=0,\ P'_n ( z_0 ) =1 \},
\end{equation}
where ${\mathcal P}_n ( {\C} )$ is the class of algebraic
polynomials of degree at most $n$, with complex coefficients. The conformal
mapping $\varphi : G \rightarrow D_{R_{0}} := \{ z : |z| < R_{0} \}$, normalized
by $\varphi (z_0)=0$ and $\varphi' (z_0)=1$, solves the same
extremal problem in the class of all analytic functions $f$ in
$G$, satisfying $f(z_0)=0$ and $f'(z_0)=1$.  Here, $R_{0}$ is the inner
conformal radius of the domain $G$ with respect to $z_0.$ Moreover,
\begin{equation} \label{5.2}
\| \varphi' -B'_n \|_{2} = \inf_{P_n \in {\mathcal P}_n (
{\C})} \{ \| \varphi' -P'_n \|_{2}\ :\ P_n(z_0) =0,\
P'_n (z_0) =1 \}
\end{equation}
(see, e.g., \cite{Ke} or \cite{Ga2}).  It is clear from (\ref{5.2}) that
\begin{equation} \label{5.3}
\lim_{n \rightarrow \infty} B'_n (z) = \varphi' (z) \quad \mbox{
and } \quad \lim_{n \rightarrow \infty} B_n (z) = \varphi
(z),\quad z \in G,
\end{equation}
because polynomials are dense in $L_2 (G)$.  A more delicate fact
of the {\it uniform} convergence of $B_n (z)$ to $\varphi (z)$ on
$\ov{G}$ was first observed by Keldysh in 1939 \cite{Ke}, for the
domains $G$ with sufficiently smooth boundaries.  He also
constructed an example of a starlike domain, bounded by a
piecewise analytic curve with one singular point, where Bieberbach
polynomials diverge.  A  considerable progress in the
area has been achieved by Mergelyan \cite{Me}, Suetin \cite{Su1},
Simonenko \cite{Si}, Andrievskii \cite{An1}-\cite{AG} and Gaier
\cite{Ga2}-\cite{Ga4}.  In particular, Andrievskii \cite{An1}
proved that the uniform convergence of Bieberbach polynomials
holds on $\ov{G}$, where $G$ is any quasidisk, and Gaier
\cite{Ga2}-\cite{Ga4} showed that the rate of this uniform
convergence is quite close to the best possible rate in uniform
polynomial approximation of the conformal mapping $\varphi$.

It is well known that a quasiconformal curve does not allow zero angles
(cusps). The first results on the uniform convergence of Bieberbach polynomials
in domains with cusps were obtained by Andrievskii \cite{An2}-\cite{An3}.
Pritsker \cite{Pr5} developed his approach to improve these results for domains
with certain interior zero angles (outward pointing cusps). The interior zero
angles seem to play a special role in this problem, as Keldysh's
counterexample, although being implicit, but gives impression that his
piecewise analytic boundary curve has an outward cusp as its only singular
point (see \cite{Ke}). We confirm this by constructing an example on the
divergence of Bieberbach polynomials at an outward pointing cusp (see Theorem
\ref{thm2} and its proof). An interesting problem arising here is to find the
critical {\em order of tangency} at this interior zero angle, separating the
convergent behavior of Bieberbach polynomials, exhibited below in Theorem
\ref{thm1}, from the divergent one for sufficiently thin cusps. This would give
a rather complete answer to the old question on the geometry of domains with
uniform convergence of Bieberbach polynomials.

Introducing the area orthonormal polynomials
$\{K_n(z)\}_{n=0}^{\infty}$, such that
\begin{equation} \label{5.4}
\int\!\!\!\int_{G} K_m (z) \overline{K_n (z)} dx dy = \left \{
\begin{array}{l} 1, \quad m=n,  \\ 0, \quad m \neq n,
\end{array} \right.
\end{equation}
one can find the following representation for Bieberbach
polynomials \cite[p. 34]{Ga1}:
\begin{equation} \label{5.5}
B_n(z)=\frac{\dis\sum_{k=0}^{n-1} \overline{K_k(z_0)} \int_{z_0}^z
K_k(t)dt}{\dis\sum_{k=0}^{n-1} |K_k(z_0)|^2},\quad n \in \N.
\end{equation}
This gives a constructive method for generating Bieberbach
polynomials via the Gram-Schmidt orthonormalization process and for
numerical approximation of the conformal mapping $\varphi$ (see
\cite{Ga1}). In addition, (\ref{5.5}) indicates the connection
with the Bergman kernel function \cite{Be}
\begin{equation} \label{5.6}
K(z,z_0)=\sum_{k=0}^{\infty} \overline{K_k(z_0)} K_k(z)=
\frac{\varphi'(z)}{\pi R_0^2}, \quad z,z_0 \in G.
\end{equation}
It is clear from (\ref{5.6}) and $\varphi' (z_0)=1$ that
\begin{equation} \label{5.7}
K(z_0,z_0)=\sum_{k=0}^{\infty} |K_k(z_0)|^2=
\frac{1}{\pi R_{0}^2}.
\end{equation}
Thus, many problems on the convergence of Bieberbach polynomials are
equivalent to those on the convergence of the integrated bilinear
series of (\ref{5.6}).

\section{Convergence and Divergence Results}

Let $\tau$ be a conformal mapping of the unit disk $D$ onto a quasidisk (cf.
\cite{LV}). We say that a Jordan arc $\gamma$ is {\it quasianalytic} if
$\gamma=\tau([-1,1])$ for such a mapping $\tau$ (quasianalytic arcs were
introduced in \cite{AG}). It is known that a quasianalytic arc is rectifiable
and quasismooth, i.e., it satisfies the following chord-arc condition, by
Lavrentiev:
$$ |\gamma(z_1,z_2)| \le M \, |z_1-z_2|, \quad \forall \, z_1, z_2 \in \gamma,$$
where $|\gamma(z_1,z_2)|$ is the length of a subarc $\gamma(z_1,z_2) \subset \gamma,$
with the endpoints $z_1,z_2,$ and $M \ge 1$ is a constant, depending only on
$\gamma.$ A Jordan curve is said to be piecewise quasianalytic, if it consists of
a finite number of quasianalytic arcs.

Suppose that $G$ is bounded by a piecewise quasianalytic curve $L=\partial G$,
with the quasianalytic arcs joining at the points $\{z_j\}_{j=1}^m \subset L$.
Two quasianalytic arcs $L_j \subset L$ and $L_{j+1} \subset L$, meeting at $z_j$,
form an $x^p$-type interior zero angle, if there exists a neighborhood of $z_j$
such that in a local coordinate system, with the origin at $z_j$, we have
$$(x,y) \in L_j \Rightarrow c_1 x^P \le y \le c_2 x^p$$
and $$(x,y) \in L_{j+1} \Rightarrow -c_2 x^p \le y \le -c_1 x^P,$$
where $P \ge p >1$ and $c_1,c_2>0.$ With these notations, our result on the
convergence of Bieberbach polynomials is stated below.

\begin{thm} \label{thm1} If $\partial G$ is piecewise quasianalytic,
with $x^p$-type interior zero angles at the joint points,
then there exist $q=q(G),\ r=r(G),\ 0<q,r<1,$ and $C=C(G)>0$ such that
\begin{equation} \label{2.1}
\| \varphi - B_n\|_{\infty} \leq C\ q^{n^r}, \quad n \in \N.
\end{equation}
\end{thm}

It is worth noting that one cannot have $r=1$ in Theorem \ref{thm1}, as
this would imply that $\varphi$ is analytic on $\ov G$, by a well known result
(see, e.g., \cite[p. 27]{Ga1}), which is obviously not the case (cf. Theorem
\ref{thm3} below).

A companion divergence result is based on Keldysh's construction in \cite{Ke},
but its geometry is made more explicit here.

\begin{thm} \label{thm2} There exists a domain with piecewise smooth boundary and
one outward pointing cusp, such that Bieberbach polynomials diverge at this cusp.
Furthermore, the boundary of this domain is analytic outside of any neighborhood
of the cusp point.
\end{thm}

One might speculate that the interior zero angle in Theorem \ref{thm2} has an
{\it exponential} order of tangency at the cusp point. It would be very
interesting to find the {\em critical order of tangency} at this interior
zero angle, separating the convergent behavior of Bieberbach
polynomials in Theorem \ref{thm1}, from the divergent one in Theorem \ref{thm2}.

\begin{rem} \label{rem0}
It is also interesting to note that the boundary of the domain in Theorem \ref{thm2} (or
in Keldysh's counterexample) cannot be piecewise analytic, as is claimed in
\cite{Ke}. In other words, the boundary arc, with endpoints meeting at the only
irregular boundary point, cannot be an analytic arc, which is the image of a segment
under a mapping, analytic in a domain containing this segment inside. Indeed,
if this arc is analytic then we can define the angle at the irregular boundary
point, by using one-sided tangents. In the case this angle is non-zero, we see that
the boundary of our domain is quasiconformal, so that the associated Bieberbach
polynomials must converge uniformly by \cite{An1}. Thus divergence is only
possible if we have a zero angle at the irregular point, which translates into
an outward pointing cusp in Keldysh's construction. However, an analytic arc
can only form an $x^p$-type zero angle, because it cannot have an arbitrarily
high order of contact, as we show in Section 4. Hence we get the uniform convergence
of Bieberbach polynomials again, by Theorem \ref{thm1}!
\end{rem}

\section{Continuity and Differentiability of a Conformal Mapping at a Cusp}

The results on the behavior of $\varphi$ at an interior zero angle, obtained in this paper,
may be of independent interest. We summarize them below for convenience of the reader.

\begin{thm} \label{thm3}
Suppose that $G$ has an $x^p$-type interior zero angle at $z \in \partial
G$, which is formed by two quasianalytic arcs.  Then there exist constants $C,c>0$ such
that
\begin{equation} \label{3.1}
|\varphi(t)-\varphi(z)| \leq C \ \exp \left(- \frac{c}{|t-z|^{p-1}}\right),
\quad t \to z,\ t \in {\ov G}.
\end{equation}

Furthermore,
\begin{equation} \label{3.2}
\lim_{t \to z \atop t \in {\ov G \setminus \{z\}}} \ \frac{\varphi^{(k)}(t)}{|t-z|^m} = 0,
\quad \forall \ k,m \in \N.
\end{equation}
\end{thm}

\begin{cor} \label{cor1}
Suppose that $\partial G$ is piecewise quasianalytic, with $x^p$-type interior
zero angles at the joint points.  Then $\varphi \in C^{\infty}({\ov G})$,
i.e.,
\begin{equation} \label{3.3}
\varphi^{(k)} \in  C({\ov G}), \quad \forall \ k \in \N.
\end{equation}
\end{cor}

\section{Proofs}

Let $C,c,c_1,c_2, \ldots$ denote positive constants, not necessarily
the same at different places. Writing $a \ole b$, we mean that
$a \leq c_1 b $ for a constant $c_1$, which doesn't depend on $a$ and $b$.
The relation $a \sim b$ indicates that $c_2 b \le a \le c_1 b$, where
$c_1,c_2$ are independent of $a$ and $b$.

Our proofs heavily rely on the distortion properties of conformal
and quasiconformal mappings, where we start with the following
lemma (see Andrievskii \cite[pp. 97-98]{ABD}).

\begin{lem} \label{lem1} Let $w =F(\zeta )$ be a $K$-quasiconformal
mapping of the plane onto itself, such that $F( \infty ) = \infty$, $\zeta_j \in {\C}$,
$w_j =F(\zeta_j)$ $(j =1,2,3)$,  and $|w_1 - w_2| \leq c_1 | w_1 -w_3|$.
Then $|\zeta_1 - \zeta_2| \leq c_2 | \zeta_1 - \zeta_3|$ and
\begin{equation}\label{4.1}
c_3 \left| \frac{w_1 - w_3}{w_1 -w_2} \right|^{1/K} \leq
\left| \frac{\zeta_1 - \zeta_3}{\zeta_1-\zeta_2} \right| \leq
c_4 \left| \frac{w_1 - w_3}{w_1 -w_2} \right|^K,
\end{equation}
where $c_j = c_j(c_1,K)$, $j =2,3,4$.
\end{lem}

We next introduce the arc $l_{q,a} \subset D$, with the endpoints at $\pm 1$,
as the union of two arcs $\{ w=x+iy :\ y=(x+1)^q/a,\ -1 < x \leq 0\}$ and $\{
w=x+iy :\ y=(1-x)^q/a,\ 0 \leq x < 1 \}, \quad q,a>1.$ If $\tau$ is the
conformal map of $D$, defining a quasianalytic arc of the boundary $L_1 \subset
\partial G$ with the cusp points $z_1=\tau(-1)$ and $z_2=\tau(1)$, then we set
$$L_{q,a} := \tau(l_{q,a}).$$ We assume here that $\tau$ is extended to a
$K$-quasiconformal homeomorphism of the complex plane onto itself,
with infinity as a fixed point. It is clear that $L_{q,a}$ is an arc connecting
$z_1$ and $z_2$ in the exterior of $G$. We give estimates for the distance
$d(\zeta,L_1)$ from $\zeta \in L_{q,a}$ to $L_1$, $q,a>1$, and also for the
distance $d(\zeta,L_{q,a})$ from $\zeta \in L_1$ to $L_{q,a}$.

\begin{lem} \label{lem2} If $\zeta \in L_{q,a}$ then
\begin{equation}\label{4.2}
\min_{j=1,2} |\zeta-z_j|^{(q-1)K^2+1} \ole d(\zeta,L_1) \ole
\min_{j=1,2} |\zeta-z_j|^{(q-1)/K^2+1}.
\end{equation}
Similarly, if $\zeta \in L_1$ then
\begin{equation}\label{4.2a}
\min_{j=1,2} |\zeta-z_j|^{(q-1)K^2+1} \ole d(\zeta,L_{q,a}) \ole
\min_{j=1,2} |\zeta-z_j|^{(q-1)/K^2+1}.
\end{equation}
\end{lem}

\begin{proof}  Set $t=\tau^{-1}(\zeta)$ and $x = {\Re} \, t$. Assume that
$x \geq 0$. Since $|t - x| < |t - 1|$, we obtain that
$$ |t-1|^{1-q} \sim \left| \frac{t-1}{t-x} \right| \ole \left| \frac{\zeta -z_2}{\zeta -\tau(x)} \right|^K \le
\left( \frac{|\zeta-z_2|}{d(\zeta,L)} \right)^K,$$
by Lemma \ref{lem1}. Therefore
$$ d(\zeta,L) \ole |t-1|^{(q-1)/K}|\zeta -z_2|.$$
Applying (\ref{4.1}) again, with $\zeta_1=z_2,\ \zeta_2=\zeta \mbox{ and } \zeta_3=z_1$, we have
\begin{equation}\label{4.3}
|\zeta-z_2|^K \ole |t-1| \ole |\zeta-z_2|^{1/K},
\end{equation}
which gives the right hand side of (\ref{4.2}), by the previous inequality.

Let $\zeta' \in L$ be such that $d(\zeta,L) = |\zeta -\zeta'|$,
and set $t' = \tau^{-1}(\zeta')$.  Assume that ${\Re}\, t' \geq 0$. Since
$| \zeta - \zeta' | \leq |\zeta -z_2|$, Lemma \ref{lem1} yields
$$\frac{|\zeta -z_2|}{d(\zeta,L)} = \left| \frac{\zeta -z_2}{\zeta-\zeta'} \right|
\ole \left| \frac{t-1}{t-t'} \right|^K \le \left| \frac{t-1}{t-x} \right|^K \sim |t-1|^{(1-q)K}.$$
Combining the above estimate with (\ref{4.3}), we also prove the left hand side of (\ref{4.2}).

The estimates in (\ref{4.2a}) are obtained by an analogous argument.
\end{proof}

\begin{rem} \label{rem1}
If we consider $\ov{l}_{q,a} :=\{ \bar{z} : z \in l_{q,a} \} \subset
D$ and set $ \ov{L}_{q,a} := \tau(\ov{l}_{q,a}),$ then (\ref{4.2}) holds for
any $\zeta \in \ov{L}_{q,a}$ as well.
\end{rem}

We now construct an analytic extension of the conformal mapping $\varphi$ into a domain
$\tilde G$ containing $G$.

\begin{lem} \label{lem3}
Suppose that $\partial G$ is piecewise quasianalytic, with $x^p$-type interior
zero angles at the joint points. Then the mapping $\varphi$ can be continued conformally into a domain
$\tilde G$, with the rectifiable boundary $\partial \tilde G$ that consists of quasismooth
arcs $L_{q,a}$, connecting the cusp points $\{z_j\}^m_{j=1}$, such that $G \subset \tilde G$ and
$\partial G \cap \partial \tilde G = \{z_j\}^m_{j=1}$.

Furthermore, there exist constants $C,c>0$ such that
\begin{equation} \label{4.4}
|\varphi(z)-\varphi(z_j)| \leq C \ \exp \left(-\frac{c}{|z-z_j|^{p-1}}\right),
\quad z \in \tilde G,
\end{equation}
where $j=1,\ldots,m.$
\end{lem}

\begin{proof}
Let $L_j$, $j =1,2, \ldots ,m,$ be the quasianalytic arc of $L = \partial G$,
connecting $z_j$ and $z_{j+1}$, which is defined by the corresponding conformal map
$\tau_j$. Denote the domain, bounded by the arcs $l_{q,a}$ and ${\ov l}_{q,a}$, by
$S_{q,a}$.  It follows from (\ref{4.2}) and Remark \ref{rem1} that
$$\tau_j (S_{q,a}) \cap \tau_k (S_{q,a}) = \emptyset \quad \mbox{ for } j
\neq k,$$
provided we choose $a>0$ to be sufficiently small and $q$ to satisfy $(q-1)/K^2+1>P$.

On defining
$$\tilde G := G \cup \left[ \bigcup_{j=1}^m \tau_j (S_{q,a}) \right],$$
we extend the conformal mapping $\varphi$ into $\tilde G$ using the standard reflection
principle:
\begin{equation}\label{4.5}
\varphi (z) :=\frac{R_{0}^2}{\overline{\varphi \left[ \tau_j \left( \overline{\tau^{-1}_j(z)}
\right) \right]} }, \quad z \in \tau_j(S_{q,a}) \backslash \ov{G},
\end{equation}
where $j=1,\ldots,m.$

We next proceed to proving (\ref{4.4}), where we use the method of moduli of
curve families (the method of extremal length). There is no loss of generality
in assuming that $z_j=0 \in \partial G$ and that $G \cap D_R(0) \subset W$ for
some $R>0$, where the wedge $W$ is defined by $W:=\{w=x+iy:|y|<c_1x^p,x>0\},$
with $p>1.$ Fix a point $a \in G$ and set $d=\min(|a|/2,R).$ Consider a point
$z \in \ov{G},$ such that $|z|<d$, and a family of curves $\Gamma$, separating
points $0$ and $z$ from the point $a$ in $G$. We need to estimate the module
of $\Gamma$, denoted by $m(\Gamma)$, from below. This is accomplished with the
help of an auxiliary family of curves $\Gamma'$, which consists of the circular
arcs $\gamma(r):=\{|w|=r\} \cap W$, with the radius $r$ varying from $c_2|z|$ to
$d$, $c_2|z|<r<d.$ Since the boundary arcs $L_j$ and $L_{j+1}$, meeting at $z_j=0$,
are quasismooth, we can choose $c_2$ such that each curve from $\Gamma'$ contains
a curve from $\Gamma.$ It follows by the comparison principle (see Theorem 4-1 in
\cite[p. 54]{Ah}) that
\begin{equation}\label{4.6}
m(\Gamma) \geq m(\Gamma').
\end{equation}
Let $\theta(r)$ be the angular measure of the arc $\gamma(r) \in \Gamma'.$ It
is known (cf. Theorem 2.6 in \cite[p. 77]{Oh}) that
\begin{equation}\label{4.7}
m(\Gamma') = \int_{c_2|z|}^d \frac{dr}{r\theta(r)}.
\end{equation}
Using a simple estimate $\theta(r) \leq 2c_1r^{p-1}, \ p>1$, we conclude by
(\ref{4.6}) and (\ref{4.7}) that
$$ m(\Gamma) \geq \int_{c_2|z|}^d \frac{dr}{2c_1r^{p}}=\frac{1}{2(p-1)c_1}
\left(\frac{1}{(c_2|z|)^{p-1}} - \frac{1}{d^{p-1}}\right).$$
Hence
\begin{equation} \label{4.8}
|\varphi(z)-\varphi(0)| \ole e^{-\pi m(\Gamma)} \ole \exp \left(-\frac{c}{|z|^{p-1}}\right),
\quad z \in {\ov G},
\end{equation}
by Theorem 1 of \cite[p. 290]{Bel} (see also \cite[p. 34]{ABD}).
In the case $z \in \tilde G \setminus {\ov G}$, we obtain from (\ref{4.5}) and (\ref{4.8})
that
$$ |\varphi(z)-\varphi(0)| \ole  \exp \left(-\frac{c}{\left| \tau_j \left(
\overline{\tau^{-1}_j(z)} \right)\right|^{p-1}}\right),$$
with a different $c$.
Applying (\ref{4.1}) with $\zeta_1=0,\ \zeta_2= \tau_j \left( \overline{\tau^{-1}_j(z)}
\right), \ \zeta_3=z$ and $w_1=-1,\ w_2=\overline{\tau^{-1}_j(z)},\
w_3=\tau^{-1}_j(z),$ we conclude that
$$ \left| \tau_j \left( \overline{\tau^{-1}_j(z)} \right) \right| \sim |z|,$$
which implies (\ref{4.4}) by the previous inequality.
\end{proof}

\begin{lem} \label{lem5}  Let $G$ be a Jordan domain, which is symmetric
in the real axis, and let $z_0 =0 \in G$.  Assume that $\xi \in \partial G$ is real and
$G \subset \{z: {\Re}\ z < \xi \}$.  If the conformal mapping $\varphi$ is not
analytic on $\ov{G}$, then
\begin{equation}\label{4.50}
\limsup_{n \rightarrow \infty} | B_n (x)| = \infty , \qquad \forall\ x > \xi .
\end{equation}
\end{lem}

\medskip\noindent {\bf Proof.}
Suppose to the contrary of (\ref{4.50}) that there exists
$x_0 > \xi$ such that the sequence $\{ B_n (x_0) \}_{n =0}^{\infty}$ is
bounded.  Then we obtain from (\ref{5.5}) and (\ref{5.7}) that the following
sequence is also bounded:
$$\sum_{k =0}^{n -1} \overline{K_k (0)} \int_0^{x_0} K_k (t)dt, \quad
n \in {\N}.$$
This implies, in turn, that
\begin{equation}\label{4.51}
\left| \overline{K_n (0)} \int_0^{x_0} K_n (t)dt \right| < C, \quad n \in {\N},
\end{equation}
for some constant $C >0$.  Note that the orthonormal polynomials $K_n (z)$
have real coefficients for any $n \in {\N}$, because $G$ is symmetric about the real axis.
Furthermore, we follow the usual convention in Gram-Schmidt orthonormalization that the
leading coefficient of $K_n (z)$ is positive for any $n \in {\N}$.
It follows that each $K_n (x)$ is real valued for real $x$, and is positive for
$x \rightarrow + \infty$. Since the zeros of $K_n (z)$
are contained in the convex hull of $\overline{G}$ (see \cite[p. 31]{StT}), we
conclude that $K_n (x)$ has no zeros for $x > \xi$ and that
\begin{equation}\label{4.52}
K_n (x)> 0, \qquad x > \xi, \qquad n \in {\N}.
\end{equation}
Using Theorem 1.1.4 of \cite[p. 4]{StT}, we obtain that
\begin{equation}\label{4.53}
\limsup_{n \rightarrow \infty} \| K_n \|_{\infty}^{1/n} \leq 1
\end{equation}
and
\begin{equation}\label{4.54}
\lim_{n \rightarrow \infty} | K_n (x)|^{1/n} = e^{g_{\Omega}(x, \infty )}
>1, \qquad x > \xi,
\end{equation}
where $g_{\Omega} (x, \infty )$ is the Green function of $\Omega := \ov{\C}
\backslash \ov{G}$ with pole at $\infty$.  Combining (\ref{4.52})-(\ref{4.54}) gives that
$$\liminf_{n \rightarrow \infty} \left| \int_0^{x_0} K_n (t)dt \right|^{1/n} >1$$
and that
\begin{equation}\label{4.55}
\limsup_{n \rightarrow \infty} | K_n (0)|^{1/n} < 1,
\end{equation}
by (\ref{4.51}).  Hence the conformal mapping $\varphi$ must have an analytic
continuation through
$\partial G$, by (\ref{4.55}) and Theorem 2.1 of \cite{PSG}, which contradicts our
assumption.\hfill\qed

\medskip\noindent {\bf Proof of Theorem \ref{thm1}.}
We use a known method, based on the extremal property of Bieberbach
polynomials (\ref{5.2}). Namely, we first estimate the quantity
$\|\varphi'-B_n'\|_2$, and then proceed to the uniform norm case, to prove
(\ref{2.1}).

Recall that the conformal mapping $\varphi$ can be continued, by Lemma
\ref{lem3}, into a larger domain $\tilde G$, whose boundary consists of
quasismooth arcs connecting the cusp points $\{z_j\}^m_{j=1} \subset
\partial G$. Let $\gamma_j$ be a subarc of $\partial \tilde G$, with the
endpoints $z_j$ and $z_{j+1}$, and let $\zeta_j \in \gamma_j$ be a fixed point,
$j=1,\ldots,m.$ Note that $\zeta_j$ divides $\gamma_j$ into  $\gamma_j^1$ and
$\gamma_j^2$, so that $\partial \tilde G = \bigcup_{j =1}^m
\bigcup_{i =1}^2 \gamma_j^i$. Since $\partial \tilde G$ is rectifiable, we have
by Cauchy's integral formula that
\begin{eqnarray*}
\varphi(z)  &=& \frac{1}{2 \pi i}  \int_{\partial \tilde G}
\frac{\varphi ( \zeta )}{\zeta -z} d \zeta = \frac{1}{2 \pi i} \sum_{j =1}^m \sum_{i =1}^2
\int_{\gamma_j^i} \frac{\varphi(\zeta )}{\zeta-z} d \zeta\\
& = & \frac{1}{2 \pi i} \sum_{j =1}^m \sum_{i =1}^2 \int_{\gamma_j^i}
\frac{\varphi(\zeta) - \varphi(z_{j +i -1})}{\zeta -z} d \zeta + \frac{1}{2 \pi
i} \sum_{j =1}^m \varphi (z_j) \log \frac{\zeta _{j -1} -z}{\zeta_j -z},
\end{eqnarray*}
for any $z \in G$, where we assume that $\zeta_0 = \zeta_m$. It follows that
\begin{eqnarray} \label{4.20}
\varphi'(z) &=& \frac{1}{2 \pi i} \sum_{j =1}^m \sum_{i =1}^2 \int_{\gamma_j^i}
\frac{\varphi(\zeta) - \varphi(z_{j +i -1})}{(\zeta -z)^2} d \zeta \\
&+& \frac{1}{2 \pi i} \sum_{j =1}^m \varphi (z_j) \left( \frac{1}{\zeta_j -z}-
\frac{1}{\zeta _{j -1}-z} \right), \quad z \in G. \nonumber
\end{eqnarray}
Observe that the second sum
$$f(z) := \frac{1}{2 \pi i} \sum_{j =1}^m \varphi (z_j) \left( \frac{1}{\zeta_j -z}-
\frac{1}{\zeta _{j -1}-z} \right)$$
represents a function, analytic on $\ov G.$ Hence there exists a sequence of
polynomials $\{p_n\}_{n=1}^{\infty}$ and a number $R>1$ such that
\begin{equation}\label{4.21}
\|f-p_n\|_{\infty} \leq C R^{-n}, \quad n \in \N,
\end{equation}
(see, e.g., Theorem 4 in \cite[p. 27]{Ga1}). Consequently, our problem of
approximating $\varphi'(z)$ by polynomials reduces to approximating functions
of the following form
$$g_{i,j}(z) := \int_{\gamma_j^i} \frac{\varphi(\zeta) - \varphi(z_{j +i -1})}
{(\zeta -z)^2} d \zeta,$$
in view of (\ref{4.20}). Furthermore, we can consider approximation in the
uniform norm and then pass to $L_2(G)$ norm, which suffices for our purposes.
We now set
\begin{equation}\label{4.22}
g(z) := g_{1,1}(z) = \int_{\gamma} \frac{\varphi(\zeta) - \varphi(z_1)}
{(\zeta -z)^2} d \zeta, \quad \gamma:=\gamma_1^1,
\end{equation}
and study the approximation of this function only, as the other functions
$g_{i,j}(z)$ are handled similarly.

Let $\Phi:\Omega \to D'$ be a conformal map of $\Omega:={\ov \C} \setminus {\ov
G}$ onto $D':=\{w:|w|>1\}$, satisfying the conditions $\Phi(\infty)=\infty$ and
$\Phi'(\infty)>0.$ Define the level curves of $\Phi$ by
$$ L_R :=\{z \in {\ov \Omega}: |\Phi(z)|=R\}, \quad R \geq 1,$$
where we set $L:=L_1=\partial \Omega =\partial G.$ Let $G_R:=Int\, L_R,\ R>1,$
be the domain bounded by $L_R.$ Clearly, if $R$ is sufficiently close to 1,
then $\gamma \cap L_R \neq \emptyset.$ Denote $\gamma':=\gamma \cap \ov{G_R}$ and
$\gamma'':=\gamma \setminus \gamma'$,  so that $\gamma''$ lies exterior to $L_R$.
Hence the function
\begin{equation}\label{4.23}
h_2(z) := \int_{\gamma''} \frac{\varphi(\zeta) - \varphi(z_1)}
{(\zeta -z)^2} d \zeta
\end{equation}
is holomorphic in $G_R$ and is well approximable by polynomials.
Namely, we obtain from Theorem 3 of \cite[p. 145]{SL} that there exists a
sequence of polynomials $\{p_n\}_{n=1}^{\infty}$ such that
\begin{equation}\label{4.24}
\|h_2-p_n\|_{\infty} \leq C \frac{n}{(r-1)^2} \max_{z \in G_r} |h_2(z)| \ r^{-n},
\quad n \in \N,
\end{equation}
where $C$ is an absolute constant and $r<R$. On choosing $R=1+2n^{-s}$ and
$r=1+n^{-s}$, with $s \in (0,1)$, we estimate
\begin{eqnarray*}
\max_{z \in G_r} |h_2(z)| &\leq& \int_{\gamma''} \frac{|\varphi(\zeta) - \varphi(z_1)|}
{|\zeta -z|^2} |d \zeta| \ole \frac{1}{\dis \min_{z \in G_r,\, \zeta \in \gamma''}
|\zeta -z|^2} \\ &\ole& \frac{1}{[d(L_R,L_r)]^2},
\end{eqnarray*}
where $d(L_R,L_r)$ is the distance between $L_R$ and $L_r.$ Note that
$$ d(L_R,L_r) \ge c(R-r)^2, $$
by a result of Loewner (see \cite[p. 61]{ABD}),
which implies
$$ d(L_R,L_r)\oge n^{-2s}. $$
We conclude that
$$ \max_{z \in G_r} |h_2(z)| \ole n^{4s},$$
and, using (\ref{4.24}),
\begin{equation}\label{4.25}
\|h_2-p_n\|_{\infty} \ole n^{1+2s+4s}(1+n^{-s})^{-n} \ole
n^{1+6s}e^{-n^{1-s}}, \quad n \in \N.
\end{equation}

Introducing a companion function
$$ h_1(z) := \int_{\gamma'} \frac{\varphi(\zeta) - \varphi(z_1)}
{(\zeta -z)^2} d \zeta, $$
so that
\begin{equation}\label{4.26}
g(z)=h_1(z)+h_2(z),
\end{equation}
we now show that $\|h_1\|_{\infty}$ is sufficiently small. Indeed, we obtain
by Lemma \ref{lem3} and \ref{lem2} that
\begin{eqnarray*}
\|h_1\|_{\infty} &\leq& \max_{z \in {\ov G}}  \int_{\gamma'} \frac{|\varphi(\zeta) - \varphi(z_1)|}
{|\zeta -z|^2} |d \zeta| \ole \int_{\gamma'} \frac{\exp\left(-c/|\zeta-z_1|^{p-1}\right)}
{[d(\zeta,L)]^2} |d \zeta| \\ &\ole& \int_{\gamma'} \frac{\exp\left(-c\, [d(\zeta,L)]^
{-\frac{p-1}{K^2(q-1)+1}}\right)}{[d(\zeta,L)]^2} |d \zeta| \\
&\ole& \max_{\zeta \in \gamma'} \frac{\exp\left(-c\, [d(\zeta,L)]^{-\frac{p-1}{K^2(q-1)+1}}\right)}
{[d(\zeta,L)]^2}.
\end{eqnarray*}
Since the function $x^{-2}\exp(-cx^{-a})$, where $a,c>0$,  is strictly
increasing on an interval $(0,x_0)$, we deduce from the previous inequality
that
\begin{equation}\label{4.27}
\|h_1\|_{\infty} \ole \frac{\exp\left(-c\, [d(\zeta_R,L)]^{-\frac{p-1}{K^2(q-1)+1}}\right)}
{[d(\zeta_R,L)]^2},
\end{equation}
where $R=1+2n^{-s}$ is sufficiently close to 1 and $\zeta_R \in L_R$. It is known that
$\Psi:=\Phi^{-1}$ is H\"{o}lder continuous on $\ov{D'}$ (see Theorem 3 in
\cite{NP}), so that
$$d(\zeta_R,L):=\min_{t \in L} |\zeta_R-t| \ole (R-1)^{\beta} \ole
n^{-s\beta},$$
for some $\beta>0.$ Hence we obtain from (\ref{4.27}) that
\begin{equation}\label{4.28}
\|h_1\|_{\infty} \ole n^{2s\beta}
\exp\left(-c\, n^{\frac{s\beta(p-1)}{K^2(q-1)+1}}\right), \quad  n \in \N.
\end{equation}
Combining (\ref{4.25}), (\ref{4.26}) and (\ref{4.28}), we have
$$ \|g-p_n\|_{\infty} \ole \exp\left(-cn^r\right), \quad n \in \N,$$
where $r \in (0,1)$ is any number satisfying
\begin{equation}\label{4.29}
r < \min\left(1-s,\frac{s\beta(p-1)}{K^2(q-1)+1}\right).
\end{equation}
Furthermore, this immediately implies that there exists a sequence of
polynomials $\{P_n(z)\}_{n=1}^{\infty}$ such that
\begin{equation}\label{4.30}
\|\varphi'-P_n\|_{2} \ole \|\varphi'-P_n\|_{\infty} \ole \exp\left(-cn^r\right), \quad n \in \N,
\end{equation}
by (\ref{4.20}). This concludes the first part of the proof, because we obtain
from (\ref{4.30}) and (\ref{5.2}) that
\begin{eqnarray*}
\|\varphi'-B'_n\|_{2} &\le& \|\varphi'-(P_{n-1}-P_{n-1}(z_0)+1)\|_{2} \\ &\ole&
\|\varphi'-P_{n-1}\|_{2} + |1-P_{n-1}(z_0)| \\ &=&
\|\varphi'-P_{n-1}\|_{2} + |\varphi'(z_0)-P_{n-1}(z_0)| \\ &\ole&
\exp\left(-c(n-1)^r\right) \ole \exp\left(-cn^r\right),
\end{eqnarray*}
where $n \ge 2.$

The second part follows from a standard argument on translating the estimate
\begin{equation}\label{4.31}
\|\varphi'-B'_n\|_{2} \ole \exp\left(-cn^r\right), \quad n \in \N,
\end{equation}
into (\ref{2.1}), using the following polynomial inequality
\begin{equation}\label{4.32}
\|Q_n\|_{\infty} \ole n^{P-1}\ \|Q'_n\|_2, \quad n \in \N,
\end{equation}
which is valid for any $Q_n(z),$ with $Q_n(z_0)=0$ (see Corollary 2 in \cite{Pr10}).
We write
$$ \varphi(z) = B_n(z) + \sum_{k=1}^{\infty} \left(B_{(k+1)n}(z) -
B_{kn}(z)\right), \quad z \in G,$$
so that
\begin{equation}\label{4.33}
\|\varphi - B_n\|_{\infty} \le \sum_{k=1}^{\infty} \|B_{(k+1)n} -
B_{kn}\|_{\infty}.
\end{equation}
We next estimate terms in the above sum, using (\ref{4.32}) and (\ref{4.31}):
\begin{eqnarray*}
\|B_{(k+1)n} - B_{kn}\|_{\infty} &\ole& \left((k+1)n\right)^{P-1} \|B'_{(k+1)n} -
B'_{kn}\|_2 \\ &\le&
\left((k+1)n\right)^{P-1} \left(\|\varphi'-B'_{kn}\|_2 +
\|B'_{(k+1)n}-\varphi'\|_2 \right) \\ &\ole&
(k+1)^{P-1}n^{P-1} \ \exp\left(-c(kn)^r\right).
\end{eqnarray*}
It follows from (\ref{4.33}) that
\begin{eqnarray*}
\|\varphi - B_n\|_{\infty} &\ole& \sum_{k=1}^{\infty} (k+1)^{P-1}n^{P-1} \
\exp\left(-c(kn)^r\right) \\ &=& n^{P-1}\ \exp\left(-cn^r\right) \sum_{k=1}^{\infty} (k+1)^{P-1} \
\exp\left(-cn^r(k^r-1)\right) \\ &\le& n^{P-1}\ \exp\left(-cn^r\right) \sum_{k=1}^{\infty} (k+1)^{P-1} \
\exp\left(-c(k^r-1)\right) \\ &\ole& n^{P-1}\ \exp\left(-cn^r\right).
\end{eqnarray*}
Since we can drop the term $n^{P-1}$ in the last estimate, by slightly decreasing
$r$, equation (\ref{2.1}) is proved.
\hfill\qed

\medskip\noindent {\bf Proof of Theorem \ref{thm2}.}
We essentially follow Keldysh's construction in \cite{Ke},
augmented with Lemma \ref{lem5}. Let $G_1 \subset \{ |z| < 1 \}$
be a symmetric in the real axis domain, which is bounded by a
piecewise analytic Jordan curve with the only corner point $\xi_1
\subset (0,1)$. Clearly, if the inner angle at $\xi_1$ is $\alpha_1
\pi$, where $\alpha_1 \in (0,1)$ is irrational, then the conformal mapping
of $G_1$ onto a disk cannot be analytic in a neighborhood of
$\xi_1$.  Therefore, $G_1$ satisfies the assumption of Lemma
\ref{lem5} and we can find a point $\xi_2 \subset ( \xi_1,1)$ and
a number $n_1 \in {\N}$ such that
$$| B_{n_1,1} ( \xi_2)| > 2,$$
where $B_{n,1} (z)$ is the $n$-th Bieberbach polynomial associated
with $G_1$. Next, we similarly construct a domain $G_2$ bounded by
a symmetric piecewise analytic curve with the only corner point at
$\xi_2$, so that $G_1 \subset G_2 \subset \{ |z| < 1 \}$ and
$$\max_{|z| \leq 1} |B_{n,1} (z) -B_{n,2} (z)| < \frac{1}{2^2},
\quad n \leq n_1.$$
This can be always achieved by taking the
boundary of $G_2$ sufficiently close to the continuum
$\overline{G_1} \cup [ \xi_1 , \xi_2]$, because the coefficients
of Bieberbach polynomials are rational functions of the moments
$\int \!\!\! \int_{G_1} z^k \ov{z}^{\ell} dxdy$, by Gram-Schmidt
orthonormalization scheme and (1.5), and are continuously
dependent on the domain.  Proceeding in this fashion, we obtain a
sequence of domains $G_1 \subset G_2 \subset \ldots \subset G_m
\subset \ldots \subset \{ |z| < 1 \}$, such that
\begin{equation}\label{4.56}
|B_{n_m ,m} ( \xi_{m +1}) > 2m, \quad m \in {\N},
\end{equation}
and
\begin{equation}\label{4.57}
\max_{|z| \leq 1} |B_{n,m} (z) -B_{n ,m+1} (z)| < \frac{1}{2^{m +1}}, \quad
n \leq n_m,
\end{equation}
where $B_{n,m} (z)$ is the $n$-th Bieberbach polynomial associated with the
domain $G_m$.  Furthermore, we can carry out this construction in such a way that
$G_m$ converges to a domain $G$, as $m \rightarrow \infty$, which is bounded by a
piecewise smooth curve symmetric in the real axis, with the only singular point
$\xi := \lim_{m \rightarrow \infty} \xi_m$.  Let $B_n(z)$ be the $n$-th
Bieberbach polynomial for $G$.  Then we have that
\begin{equation}\label{4.58}
\lim_{m \rightarrow \infty} B_{n,m} (z) = B_n (z), \quad n \in {\N},
\end{equation}
where the convergence is uniform on compact subsets of ${\C}$ for each
fixed $n \in {\N}$.  It follows from (\ref{4.57}) and (\ref{4.58}) that
$$\max_{|z| \leq 1} |B_{n_m ,m} (z) -B_{n_m} (z)| < \frac{1}{2^m}, \quad m
\in {\N},$$
which implies that
$$|B_{n_m} ( \xi_{m +1})| > m, \quad m \in {\N},$$
by (\ref{4.56}).  Hence
$$\limsup_{n \rightarrow \infty} \| B_n \|_{\infty} = \infty .$$
But this is impossible if $\partial G$ has a non-zero angle at $\xi$, as
Bieberbach polynomials converge uniformly on $\ov{G}$, for domains with quasiconformal
boundary (see \cite{An1}).  Thus, we are forced to conclude that $\partial G$ has an
outward pointing cusp at $\xi$, according to our construction.

To show that $\partial G$ is analytic outside of any neighborhood of $\xi$, we
specify the construction of the domains $G_m$ as follows. The piecewise
analytic boundary of $G_m$, with a corner at $\xi_m$, is defined by $\partial G_m =
\tau_m([-1,1]),\ \tau_m(-1)=\tau_m(1)=\xi_m,$ for a mapping $\tau_m$ analytic in
$\{w:|w|<1+\eps_m\},$ where $\eps_m \searrow 0,$ as $m \to \infty.$ Clearly, each
$\tau_m$ is bounded in the unit disk $D$, for any $m \in \N.$ Therefore we can
find a subsequence $\tau_{m_k}$ that converge locally uniformly in $D$ to an
analytic mapping $\tau$, by a normal families argument. It follows from our
geometric construction that $\partial G = \tau([-1,1]),\ \tau(-1)=\tau(1)=\xi.$
\hfill\qed

\medskip\noindent {\bf Proof of Remark \ref{rem0}.}
We need to show that an analytic arc, which is different from a segment of the real
axis, can only have a finite order of contact
with the real axis. It is sufficient to consider the case of an analytic mapping
$\tau : [0,1] \to \gamma$, defining the arc $\gamma$, such that
$$ \tau(w)=\sum_{n=1}^{\infty} c_n w^n,$$
where this series converges in a neighborhood of $w=0.$ Suppose that
$a_n=\Re\,c_n$ and $b_n=\Im\,c_n$. Hence
$$ x(t):=\Re\, \tau(t)=\sum_{n=1}^{\infty} a_n t^n \quad \mbox{ and } \quad
y(t):=\Im \,\tau(t)=\sum_{n=1}^{\infty} b_n t^n,$$
for $t \in [0,\eps).$ If the arc $\gamma$ has higher order of tangency than at
any $x^p$-type zero angle at $\tau(0)=0$, then
$$\lim_{t \to 0^+} \frac{y(t)}{(x(t))^p} = 0, \quad \forall\, p \in \N.$$
Since $a_n \neq 0$ for some $n \in \N$, we obtain that
$$\lim_{t \to 0^+} \frac{y(t)}{t^p} = 0, \quad \forall\, p \in \N.$$
Consequently, $b_n=0,\ \forall\, n \in \N,$ and $y(t)=\Im\, \tau(t) \equiv 0$ for
$t \in [0,\eps).$ We are forced to conclude that $\gamma$ is a subset of the real axis
in a neighborhood of $w=0$, which is an obvious contradiction.
\hfill\qed

\medskip\noindent {\bf Proof of Theorem \ref{thm3}.}
Note that (\ref{3.1}) follows directly from (\ref{4.4}) in Lemma \ref{lem3}. To
prove (\ref{3.2}), we consider the analytic continuation of the conformal
mapping $\varphi$, constructed in Lemma \ref{lem3}. Thus $\varphi$ is
analytic in a larger domain $\tilde G$, such that
\begin{equation} \label{4.62}
|\varphi(t)-\varphi(z)| \ole \exp \left(- \frac{c}{|t-z|^{p-1}}\right),\quad t \in \tilde
G,\ |t-z|<d,
\end{equation}
where $d>0$ is sufficiently small. Moreover, (\ref{4.2a}) gives the
following estimate for the distance from $t$ to $\partial \tilde G$:
\begin{equation} \label{4.63}
d(t,\partial \tilde G) \oge |t-z|^a,\quad t \in \ov G,\ 0<|t-z|<d,
\end{equation}
where $a > 1.$ Letting $t \in \ov G,\ 0<|t-z|<d$, we write
$$ \varphi^{(k)}(t)=\frac{k!}{2\pi i} \int_{|w-t|=r} \frac{\varphi(w)\,
dw}{(w-t)^{k+1}}, \quad k \in \N,$$
and estimate by (\ref{4.62}):
\begin{eqnarray*}
|\varphi^{(k)}(t)| &=& \left| \frac{k!}{2\pi i} \int_{|w-t|=r} \frac{\varphi(w)\,
dw}{(w-t)^{k+1}} - \frac{k!}{2\pi i} \int_{|w-t|=r} \frac{\varphi(z)\,
dw}{(w-t)^{k+1}} \right| \\ &\le& \frac{k!}{2\pi} \int_{|w-t|=r} \frac{|\varphi(w)-
\varphi(z)|\, |dw|}{r^{k+1}} \ole \frac{k!}{r^k} \exp \left(- \frac{c}{|w-z|^{p-1}}\right)
\\ &\le& \frac{k!}{r^k} \exp \left(- \frac{c}{(r+|t-z|)^{p-1}}\right).
\end{eqnarray*}
Observe that we can use any $r<d(t,\partial \tilde G)$, which implies with $r \sim |t-z|^a$
that
$$ |\varphi^{(k)}(t)| \ole |t-z|^{-ka} \exp \left(- c/|t-z|^{p-1}\right),
\quad t \in \ov G,\ 0<|t-z|<d.$$
Hence
$$ \lim_{t \to z \atop  t \in \ov G} \frac{|\varphi^{(k)}(t)|}{|t-z|^m} = 0,
\quad k,m \in \N, $$
so that (\ref{3.2}) follows.
\hfill\qed

\medskip\noindent {\bf Proof of Corollary \ref{cor1}.}
It is clear that $\varphi$ is analytic in a neighborhood of every $z \in
\partial G$, which is not a cusp point, by the analytic continuation construction of Lemma
\ref{lem3}. On the other hand, if $ z \in \partial G$ is at a cusp, then we let
$\varphi^{(k)}(z)=0, \ \forall\, k \in \N,$ so that $\varphi^{(k)}$ is
continuous at $z$ by (\ref{3.2}) of Theorem \ref{thm3}. It follows that
$\varphi^{(k)} \in C(\ov G), \ \forall\, k \in \N,$ i.e., (\ref{3.3}) holds
true.
\hfill\qed


\end{document}